\documentclass[11pt]{amsart}
\usepackage{amsthm}
\usepackage{amsopn}
\usepackage{mathrsfs}
\usepackage{graphicx}
\usepackage{mathrsfs} 
\usepackage{amsmath}

\usepackage[latin1]{inputenc}
\usepackage{amssymb} % \usepackage{amsfonts}
 
\def \yy#1\par{\centerline{\bf #1}}

\def\|#1|{{\rm univ}(#1)}

\theoremstyle{plain}
\newtheorem{thm}{Theorem}
\newtheorem{theorem}[thm]{Theorem}
\def\xx#1 {\newtheorem{#1}[thm]{#1}}

\xx Theorem
\xx Corollary
\xx Lemma

\theoremstyle{definition}

\xx Definition
\xx Example
\xx Examples
\xx Notation 
\xx Proposition
\xx Remark
\xx Convention

\renewcommand{\L}{{\mathscr L}}
\newcommand{\LL}{{\mathbb L}} 

\newcommand{\U}{{\mathscr U}} 
\newcommand{\N}{{\mathbb N}}

\renewcommand{\P}{{\mathscr P}}

\newcommand{\NN}{{\mathscr N}}

\newcommand{\K}{{\mathscr K}}
\newcommand{\Lfin}{{\LL_{\on \NN\!,\,\rm fin}}}
\newcommand{\lf}{_{\on \NN\!,\,\rm locfin}}
\newcommand{\on}{{\upharpoonright}}
\newcommand{\tK}{{\LL_{\on\NN\! ,\,w}}} 
\newcommand{\LLn}{\LL_{\on \NN}}
\newcommand{\atopp}[2]{\genfrac{}{}{0pt}{}{#1}{#2}}
\newcommand{\zo}{\{0,1\}}

\author{Martin Goldstern}
\address{Algebra/DMG, TU Wien, Wiedner Hauptstr 8-10/104, 1040 Wien, Austria}
\email{goldstern@tuwien.ac.at}
\urladdr{http://www.tuwien.ac.at/goldstern/}
\title{The typical countable algebra}
\thanks{Partially supported by the Austrian Science Foundation (FWF), grant 
P17627-N12.  \\
A preprint of this paper is available at \tt arXiv.org.}

\usepackage{aaa73}
\AAAsubmitted{September 26, 2007; December 21, 2007; January 8, 2008}
\AAApage{1}
\begin{document}

\begin{abstract}
We argue that it makes sense to talk about ``typical'' properties
of lattices, and then show that  there is, up to isomorphism,
a unique countable lattice $L^*$ (the Fra{\"\i}ss{\'e} limit
of the class of finite lattices)
that has all ``typical'' properties. 
Among these properties are:  $L^*$ is simple and locally finite,
every order preserving function can be interpolated by a lattice 
polynomial, and every finite lattice or countable 
locally finite lattice embeds into $L^*$. 

The same arguments apply to other classes of algebras assuming they have a 
Fra{\"\i}ss\'e limit and satisfy
the finite embeddability property.
\end{abstract}

\maketitle

% \medskip \hrule\medskip 
% \begin{center}{\Large\bf DRAFT}\end{center}
% \bigskip \hrule\medskip 

% \tableofcontents

 % \newpage
\section{What is a ``typical'' algebra?} 
\subsection{Lattices} 
% We first  consider countably infinite lattices. 
  What is a ``typical'' property of a countably infinite 
  lattice?   It seems clear
that the typical lattice will not be a chain, and will in fact have 
arbitrarily large finite antichains.  Similarly, a ``typical'' 
lattice  will contain arbitrarily long finite chains.   It can also 
be argued that a typical lattice should not be distributive, since
distributivity is a very special property for a lattice to have. 

But why do we consider distributivity as a more special property
than nondistributivity?   One could argue that distributivity is a
 positive property, and that all ``positive'' properties 
are special;  however, 
 distributivity can  also be seen as a negative property:  neither
$N_5$ nor $M_3$ embeds (see \cite[II.1, Theorem 1]{Gratzer:1998}).

We will define a property as ``special'' if only few countable lattices 
possess this property.  Hence we need a means of measuring how 
``large'' an infinite family of countable lattices is; 
as in \cite{most} we will use the topological notion of ``first category''
or ``meagerness'' for this purpose.

\subsection{General algebras}
Our considerations will apply to a large class of algebras.   Let 
$\LL$ be the set of countably infinite algebras of some variety. 
We want to describe the ``typical'' 
member of~$\LL$.

\begin{Definition}\label{finite}
Let $\LL$ be a class of countable algebras. 
A finite $\LL$-algebra is any finite subalgebra of an algebra in $\LL$.
\end{Definition}
(Note that there may be finite algebras in the variety generated by $\LL$
which are not finite $\LL$-algebras, e.g., if $\LL$ is the class of
all countable Boolean algebras.)

% \begin{Definition}\label{adequate}
% We call a class  $\LL$ of algebras \emph{adequate} if
% it has the following properties: 
  % % \footnote{What we really need is that the finite 
  % % elements of~$\LL$ are closed under subalgebras, and that an 
  % % increasing union of finite algebras in~$\LL$ is again an algebra 
  % %in $\LL$} 
% $\LL$ is  closed under subalgebras, $\LL$ contains
% nontrivial finite algebras, 
% the finite algebras of~$\LL$ are 
% closed under finite products, and every increasing union of finite 
% algebras from $\LL$ is again in $\LL$.
% \end{Definition}

On the set 
of all countable algebras in~$\LL$ 
whose universe is equal to a fixed set $\N= \{0,1,2,\ldots\}$ we will define a
natural topology induced by a complete metric (Proposition \ref{topology}).   
By Baire's theorem, no nonempty open set is meager (=of first category).
We say that ``\emph{almost 
all} countable algebras in~$\LL$ have a certain property'' if 
the set of all $A$ in our space  having this property 
is \emph{residual} (i.e., has meager complement). 

\begin{Theorem}
Assume that the finite $\LL$-algebras  are a Fra{\"\i}ss\'e class (i.e., 
they have the amalgamation property (Definition~\ref{AP})
and the joint embedding property (Definition~\ref{JEP}, see also 
Definition~\ref{F.class})).
Assume further that $\LL$ has the finite embeddability
property (every partial finite $\LL$-algebra embeds into a
finite $\LL$-algebra, see Definition~\ref{finite.embeddability}).   

Let $L^*$ be the Fra{\"\i}ss\'e limit (Theorem~\ref{F.limit})
of the finite $\LL$-algebras.

Then almost all algebras in~$\LL$ are isomorphic to $L^*$.
% In other words:  For any 
% property $X$ of countable lattices which is invariant under isomorphisms:
% \begin{quote}
% If $X$ is ``typical'' (i.e., the set  of lattices with $X$
% is residual),\\
% then $L^*$ also has property $X$. 
% \end{quote}
\end{Theorem}
\begin{Corollary}  Let $L^*$ be as above.  Let $X$  be any property of algebras
that is invariant under isomorphisms.  Then there are two cases: 
\begin{itemize}
\item [(a)]  Almost all countable algebras in $\LL$ (including $L^*$) have
property~$X$.  We call such a property ``typical''.
\item [(b)]   Otherwise;  then $L^*$ does not have property $X$, so almost no
countable algebra in $\LL$ has property~$X$. 
\end{itemize} Hence $L^*$ has all ``typical'' properties.
\end{Corollary}

\begin{Remark} A variant of the above theorem for relational structures 
is well known
under the keyword \emph{ubiquitous}.  
\cite{Bankston+Ruitenburg:1990} credit the 
idea to a 1984 seminar talk by Peter Cameron.   This notion and its variants have not only appeared in model theory but also in theoretical computer
science, see \cite{Droste+Kuske:2007}. 

Algebras (and partial algebras)
can of course be seen as relational structures
(by replacing each $n$-ary function by the corresponding $(n+1)$-ary
relation, see \cite[7.7]{Bankston+Ruitenburg:1990}).

Proofs are often easier if we deal with purely relational structures only: 
on an abstract level,  the set of all countable structures for a
relational language, equipped with the natural (Tychonoff) topology, 
is a \emph{compact} metric space, whereas the space of all countable
algebras (see Proposition~\ref{topology}) is not compact.  On a technical 
level, the Fra{\"\i}ss\'e  limit (see  Theorem~\ref{F.limit}) 
of a class of relational structures can be constructed as 
an increasing union of successive
one-point extensions; in the algebraic setting
this is not possible, unless we are ready to deal with 
partial algebras.   The finite embedding property
(see Definition~\ref{finite.embedding}) seems to be crucial 
for algebras, whereas it is irrelevant for relational structures, as
relational structures are trivially locally finite. 

The purpose 
of this note is to  de-emphasise the role of compactness 
and to underline
the role of the finite embeddability property in the algebraic context. 
\end{Remark} 

In Section \ref{notions} we explain (or recall) the definitions of a
Fra{\"\i}ss\'e class and Fra{\"\i}ss\'e limit, and the finite embeddability
property.  Both notions are exemplified for the case of lattices. 
Rather than giving full proofs, we sketch proofs and mention the 
main ideas; the details can be found in textbooks  (\cite{Hodges:1997}
for model theoretic notions such as the Fra{\"\i}ss\'e limit, 
\cite{Gratzer:1998} for lattice theory). 
 
 We prove the theorem in Section \ref{proofs}.  
Section \ref{examples}  describes the example of lattices and bounded lattices. 
 
\section{Basic notions} \label{notions}

\subsection{Ideals} 

Let $U$  be any set, $I \subseteq \P(U)$ a proper ideal.
 % , which is defined ``naturally''.  

% {\small (What is ``natural''?  You will know it when you see it.)}
\begin{Definition}
 We say that ``almost all'' $u\in U$ have a property $X$ iff
the set
$$ \{u \in U : \mbox{$u$ does not have property  $X$}\}$$
is in the ideal~$I$. 
\end{Definition}

For definiteness we should say ``$I$-almost all'' rather than just
``almost all''.  However, if the ideal $I$ has a sufficiently 
natural definition, we may omit  mentioning it.  (The ideal we will 
use is the ideal of meager sets on a certain product space, see
Proposition~\ref{topology}.)

\begin{Example}

 Let $U  = [0,1]$. Let $I$  be the ideal of sets of
Lebesgue measure zero.  
% {\small This is a very natural ideal to consider.}

Then, for example, 
\begin{quote}
``almost all $x\in [0,1]$ are irrational'',
\end{quote}
{\small since the set of rational numbers has measure zero; a randomly picked number
is unlikely to be rational.}
\end{Example}

\subsection{Setup}
Throughout this paper we fix a similarity type (also called signature) 
of algebras.  For notational simplicity only 
we will assume that this type is (2,2),
i.e., we will consider only algebras with two binary operations.  We call 
these two operations $\vee$ and $\wedge$.  
Later (see Section~\ref{metric.space}) 
we will restrict our attention to algebras whose universe is a subset of
a fixed countable set~$\N$. 

We will also fix a 
class $\LL$   of countable algebras (typically the countable algebras 
of a variety).

\subsection{Partial algebras  and the finite embeddability property}

\subsubsection{Weak and relative subalgebras}

An algebra $A$ is a  set $\|A|$ (called the universe of~$A$)
equipped with  two binary operations $\vee_A$ and $\wedge_A$. 
(We allow the universe of an 
algebra to be empty.) 

A partial algebra is a set equipped with two (possibly) partial operations. 
The following definition is from \cite{Gratzer:1979}. 

\begin{Definition}   Let $A$ be an algebra (or a partial algebra),
and let $S \subseteq \|A|$ be a
subset of the universe of~$A$.  We write $A\on S$ for the partial 
algebra that $A$ induces on~$S$, and we call $A\on S$ the 
\emph{relative subalgebra}
of~$A$ determined by~$S$.     Thus, a relative subalgebra of~$A$ is 
any partial algebra $B$ whose universe $\|B|$ 
is a subset of the universe of~$A$ such that 
$$a\wedge_A b=c \ \Leftrightarrow \ a\wedge _{A\on S} b = c \quad 
\mbox{(and similarly for~$\vee$)}$$
is satisfied for all $a,b,c\in \|B|$.

A \emph{weak subalgebra} of~$A$ is any partial algebra $B$ whose universe is 
a subset of~$\|A|$ and whose operations, whenever defined, agree with the 
corresponding operations on~$A$.   E.g., whenever $a,b,c\in \|B|$ and $a\wedge_B
b= c$, then also $a \wedge_A b $ is defined and equal to~$c$. 
(In contrast to the case of relative subalgebras, we allow
here that $a\wedge_B b$  is undefined even when $a\wedge_A b $ is defined and 
an element of~$ \|B|$.)

We write $B\le_w A $ if $B$ is a weak subalgebra of~$A$,
and  $B\le_r A $ if $B$ is a relative subalgebra of~$A$. 

If $B$ is a total algebra, then $B\le_w A$ is equivalent to~$B\le_r A$; 
for total algebras $A,B$ we abbreviate $B\le_w A$ to~$B\le A$.

If  $
\delta: A\to B$ is  a total isomorphism between partial algebras 
$A$ and $B$, and $C$ a partial algebra with $ \|B| \subseteq \|C| $, 
 then
we write $\delta: A \hookrightarrow_w C$ iff $B  \le_w C$, and
we say that $\delta$ embeds  $A$ weakly into~$C$.
The notation $A\hookrightarrow_w C$ means that there exists $\delta$ with 
$\delta:A\hookrightarrow_w C$.   The relations
$\delta:A\hookrightarrow_r C$  and
$A\hookrightarrow_r C$ are
defined similarly. 
\end{Definition}

\begin{Remark} From the algebraic and model-theoretic point of 
view, the notion of a relative subalgebra seems to be more natural than the 
notion of a weak subalgebra.  However, from the topological point of
view the weak subalgebras are more practical, since they correspond 
directly to a natural basis of the relevant topology (see
Proposition~\ref{BASIS}).   Lemma~\ref{lem.equiv} shows that the
distinction is irrelevant for our purposes. 

\end{Remark}

\subsubsection{Finite embeddability}

\begin{Lemma}\label{lem.equiv}
Let $\LL$ be a class of algebras which is closed under finite products and
under isomorphic images. 
% (typically a variety). 
% 
% CLOSED UNDER SUB? 
% 
Then the following are equivalent: 
\begin{itemize}
\item[(r,r)] 
Every finite \emph{relative} subalgebra of some algebra in~$\LL$
is a  \emph{relative} subalgebra of some finite algebra in~$\LL$. 
\item [(r,w)] 
Every finite \emph{relative} subalgebra of some algebra in~$\LL$
is a  \emph{weak} subalgebra of some finite algebra in~$\LL$. 
\item [(w,w)] 
Every finite \emph{weak} subalgebra of some algebra in~$\LL$
is a  \emph{weak} subalgebra of some finite algebra in~$\LL$. 
\end{itemize} 
\end{Lemma}
\begin{proof}
Trivially, the property (r,r) implies (r,w), as does the property (w,w). 

It is also easy to see that (r,w) implies (w,w):  Let $A \le_w B \in \LL$.
Then we can find an algebra $A'$ with the same universe as $A$ which 
satisfies $A \le_w A' \le_r B$.  From (r,w) we get a finite 
algebra $E\in \LL$ with $A\le_w A' \le_r E$, which  implies $A \le_w E$. 

We now prove  that 
(w,w) implies (r,r).
Let $A \le_r B\in \LL$.    As $\LL$ is closed under 
isomorphism, it is enough to find a finite algebra $E$ with $A\hookrightarrow_r E\in \LL$. \\
Let $$ I:= \{ a\wedge a', a\vee a' : a,a'\in \|A| \} \setminus \|A|. $$
Note that $I \subseteq \|B|\setminus \|A|$. 
We may assume that $I\not=\emptyset$.  [If $I=\emptyset$, then $A$ is a
total algebra, so $A \le_w E$ trivially implies $A\le_r E$ for any finite 
$E\in \LL$.]\\
For each $i\in I$
we define a partial algebra  $A_i $ with universe $\|A| \cup \{i\}$ 
such that $A \le_r A_i \le_w B$; the only operations 
which are defined in $A_i$ but not in $A$ are those operations 
which map pairs of elements of~$A$ to~$i$.   For each $i\in I$ we can 
(using (w,w))
find a finite algebra $E_i\in \LL$ such that $A_i \le_w E_i$.  Let
$E:= \prod_{i\in I} E_i$.

Let $\delta: \|A|\to \|E| $ be the natural embedding, defined
by $\delta(a) = (a,\ldots, a)$. Clearly $\delta:A \hookrightarrow_w E$,
but we claim 
that even $\delta:A\hookrightarrow_r E$ holds. 

So let $a,a'\in \|A|$ and assume that 
(wlog) $a\wedge_A  a'$ is undefined.   We have to check that 
$\delta(a)\wedge_E \delta(a')\notin \delta( \|A|)$. 

Let $i:= a \wedge_B a'$.  
Then $i \in I$.  
The $i$-th component of~$\delta(a)\wedge_E \delta(a')\in \|E|$ 
is equal to $ a \wedge_{E_i} a' = i$,
hence $\delta(a)\wedge_E \delta(a') $ is not in the 
range of~$\delta$. 
\end{proof}

\begin{Definition}\label{finite.embeddability}\label{finite.embedding}
Let $\LL$ be a class of algebras which is closed under isomorphic images
and under finite products. 
Following  \cite{Lindner-Evans:1977} we say that 
 $\LL$ has the \emph{finite
embeddability property}  (f.e.p.)\ iff every finite
\emph{weak} subalgebra of some algebra in $\LL$ is 
a \emph{weak} subalgebra of some finite algebra in $\LL$.   (Or equivalently, 
iff  every finite
\emph{relative} subalgebra of some algebra in $\LL$  
is a  \emph{relative} subalgebra of some finite algebra in~$\LL$.)
\end{Definition}

\begin{Examples} 
Every locally finite variety has the f.e.p.  
We will see below that also the variety of lattices has the finite 
embeddability property.

I am grateful to P\'eter P.~P\'alfy for pointing out the following example:
\begin{quotation}
Let $G$ be the finitely presented group $\langle\, a,b \mid b^2a =
ab^3\,\rangle$.  This group is non-Hopfian
\cite{Baumslag-Solitar:1962}, hence not residually finite
\cite[6.1.11]{Robinson:1996}.  So there is an element $g\in \|G|$
which is contained in every normal subgroup of finite index.  Let $B$
be a finite partial subgroup of~$G$ which is generated by $a$ and $b$ and 
contains $g$,  such that $b^2a=ab^3$ can be computed within~$B$.  If
$\delta:B\to E$ is a homomorphism onto a finite group~$E$, then $E$
would be generated by $\delta(a)$ and $\delta(b)$ and satisfy
$\delta(b)^2\delta(a)=\delta(a)\delta(b)^3$, so $E$ would have to be a
homomorphic image of~$G$, and hence satisfy $\delta(g)=1$.  \\ This
shows that the variety of groups does not have the f.e.p.
\end{quotation}

The following  is  an ad hoc
example of an
(admittedly artificial) 
variety where the failure of the f.e.p.\ is more obvious: 
\begin{quotation}
We have  a  binary operation
$*$ and three unary operations $p$, $q$ and~$F$.\\
The equations of the variety say that
on the range of $F$,
$*$ is a bijection with inverses $p$ and $q$:
          $$ p((Fx)*(Fy)) \approx Fx \approx  q((Fy)*(Fx)), \quad
	    p( Fz)*q(Fz) \approx Fz.$$
Then all finite algebras of the variety satisfy $Fx\approx Fy$. \\
So the  2-element partial algebra $\{ a,b\}$ with $Fa=a$, $Fb=b$ (and 
$p,q,{*}$ undefined) 
is not a relative partial
subalgebra of any finite algebra of the variety, although it
is a relative partial subalgebra of  some infinite algebra of the variety. 
\end{quotation}
% \begin{quotation}
% We have two binary operations 
% $\cdot$ and $*$, two unary operations $p$ and $q$, and two constants 
% $0$ and $c$.\\
% The equations of the variety say that
% \begin{itemize}
% \item  $\cdot$ is a semilattice operation with  least element $0 $
% ($0\cdot x=0$).  We write $xy$ for $x\cdot y$. 
% \item  $*$ is a bijection with inverses $p$ and $q$ as far as elements 
% below $c$ are concerned, i.e., 
	  % $ p((cx)*(cy)) = cx =  q((cy)*(cx))$,
  % $p( cz)*q(cz) = cz$. 
% \end{itemize} 
% Then all finite algebras of the variety satisfy the equation $c=0$. \\
% So the finite partial algebra $\{0\}$ in which $c$ is undefined  (but $\cdot$ 
% and $*$ are total) is not
% a relative subalgebra of any finite algebra of the variety, although it
% is a relative subalgebra of  some infinite algebra of the variety. 
% \end{quotation}
\end{Examples}

% \newpage
\subsubsection{Funayama's theorem}

\cite{Funayama:1953} (see also \cite[I.5, Theorem 20]{Gratzer:1998})
characterized the \emph{partial lattices} (i.e., the partial 
algebras which are relative subalgebras of lattices) as those structures 
on which there is a well-behaved notion of ``ideal'' and shows that 
on each partial lattice $P$ the map $x\mapsto (x]$ (which 
sends each $x$ to the principal ideal generated by $x$) maps
$P$  isomorphically onto a relative sublattice of 
the lattice of ideals on~$P$.    

In particular (see \cite[I.5, Lemma 21]{Gratzer:1998}): 
\begin{Proposition}\label{21}
Let $(P,\le)$ be a partial order, and define
partial binary operations $\vee$ and $\wedge$ by 
$$ a \vee b := \sup(a,b) \qquad a \wedge b:= \inf(a,b) \qquad \mbox{(whenever this is well-defined).}$$
Then the partial algebra $(P, \wedge, \vee)$ is a partial lattice.
\end{Proposition}

If $P$ is finite then also the ideal lattice over $P$ is finite, so as 
 a corollary to Funayama's theorem we get: 

\begin{Proposition} The variety of lattices has the f.e.p.
That is:  Whenever $P = (\|P|, \vee_P, \wedge_P)$ is a
finite relative subalgebra
of a lattice, then $P$ is a relative subalgebra of a 
\emph{finite} lattice. 
\end{Proposition}

\subsection{Fra{\"\i}ss\'e classes}
\subsubsection{Amalgamation and joint embedding}
We will need the following notation: 
If $A,B, B'$ are algebras
with $A\le B$ and $A\le B'$, we  write
$B\simeq_A B'$ (``$B$ is isomorphic to~$B'$ over~$A$'')
iff there exists an isomorphism $g:B\to B'$ which is 
the identity map on~$A$.

\begin{figure}[h]
\begin{picture}(0,0)%
\includegraphics{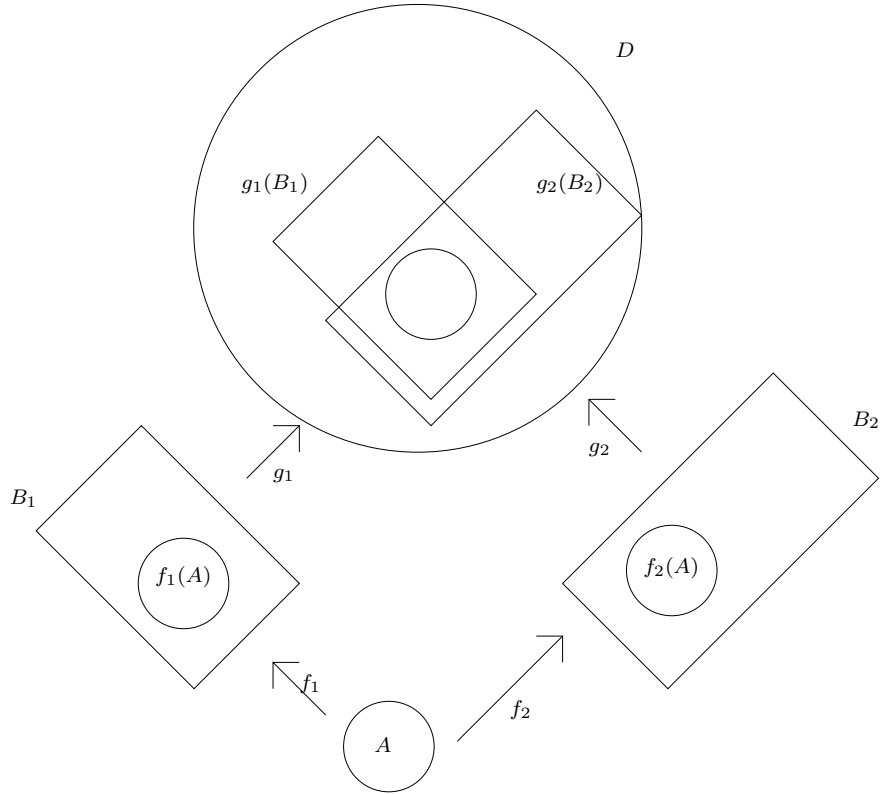}%
\end{picture}%
\setlength{\unitlength}{2901sp}%
\begingroup\makeatletter\ifx\SetFigFont\undefined%
\gdef\SetFigFont#1#2#3#4#5{%
  \reset@font\fontsize{#1}{#2pt}%
  \fontfamily{#3}\fontseries{#4}\fontshape{#5}%
  \selectfont}%
\fi\endgroup%
\begin{picture}(7452,6751)(2461,-6351)
\put(7890,-4470){\makebox(0,0)[lb]{\smash{{\SetFigFont{8}{9.6}{\rmdefault}{\mddefault}{\updefault}{$f_2(A)$}%
}}}}
\put(3720,-4550){\makebox(0,0)[lb]{\smash{{\SetFigFont{8}{9.6}{\rmdefault}{\mddefault}{\updefault}{$f_1(A)$}%
}}}}
\put(2476,-3886){\makebox(0,0)[lb]{\smash{{\SetFigFont{8}{9.6}{\rmdefault}{\mddefault}{\updefault}{$B_1$}%
}}}}
\put(9676,-3211){\makebox(0,0)[lb]{\smash{{\SetFigFont{8}{9.6}{\rmdefault}{\mddefault}{\updefault}{$B_2$}%
}}}}
\put(6751,-5686){\makebox(0,0)[lb]{\smash{{\SetFigFont{8}{9.6}{\rmdefault}{\mddefault}{\updefault}{$f_2$}%
}}}}
\put(4951,-5461){\makebox(0,0)[lb]{\smash{{\SetFigFont{8}{9.6}{\rmdefault}{\mddefault}{\updefault}{$f_1$}%
}}}}
\put(4726,-3661){\makebox(0,0)[lb]{\smash{{\SetFigFont{8}{9.6}{\rmdefault}{\mddefault}{\updefault}{$g_1$}%
}}}}
\put(7426,-3436){\makebox(0,0)[lb]{\smash{{\SetFigFont{8}{9.6}{\rmdefault}{\mddefault}{\updefault}{$g_2$}%
}}}}
\put(7651,-61){\makebox(0,0)[lb]{\smash{{\SetFigFont{8}{9.6}{\rmdefault}{\mddefault}{\updefault}{$D$}%
}}}}
\put(6976,-1186){\makebox(0,0)[lb]{\smash{{\SetFigFont{8}{9.6}{\rmdefault}{\mddefault}{\updefault}{$g_2(B_2)$}%
}}}}
\put(4456,-1186){\makebox(0,0)[lb]{\smash{{\SetFigFont{8}{9.6}{\rmdefault}{\mddefault}{\updefault}{$g_1(B_1)$}%
}}}}
\put(5591,-6006){\makebox(0,0)[lb]{\smash{{\SetFigFont{8}{9.6}{\rmdefault}{\mddefault}{\updefault}{$A$}%
}}}}
\end{picture}%
\caption[fig.amalg]{Amalgamation of~$B_1$ and~$B_2$ over~$A$}\label{fig.amalg}
\end{figure}

\begin{Definition}\label{AP}
Let $ A,  B_1,  B_2, D$ be algebras 
and let $f_1: A\to  B_1$ and 
$f_2:  A\to  B_2$ be 1-1 homomorphisms.

We say that $D$ is an amalgamation of~$B_1$ and $B_2$ over $A$ if there 
are 1-1 homomorphisms $g_1:A_1\to D$, $g_2:A_2\to D$ such that 
$g_1\circ f_1 = g_2\circ f_2$.  (See Figure~\ref{fig.amalg}.)

(More precisely we say that $(D,g_1,g_2)$ is an amalgamation 
of~$B_1$ and~$B_2$ over~$A$, or over $f_1,f_2$.
Intuitively it means that $B_1$ and
$B_2$ are glued together, identifying $f_1(A) $ with $f_2(A)$.)
\end{Definition}

\goodbreak
The following fact was noted in \cite{Jonsson:1956}: 
\begin{Proposition} \label{lat.amalg}
The family of all  lattices has the \emph{amalgamation property}, that is: 
\\
Whenever $f_1:A \to B_1$ and $f_2: A\to B_2$ are 1-1 homomorphisms between lattices,  then there exists a lattice $D$ which amalgamates $B_1$ and
$B_2$ over $A$. 
\end{Proposition}
This can be proved using 
 Funayama's characterization: We may assume that $f_1$ and 
$f_2$ are inclusion maps, and that $A = B_1\cap B_2$.  The set $ \|B_1|\cup \|B_2|$
is naturally partially ordered by the transitive closure of the union 
of the orders $\le_{B_1}$ and $\le_{B_2}$; by Proposition~\ref{21},
the partial 
operations $\sup(x,y)$ and $\inf(x,y)$ make this poset into a partial 
lattice $B_1\cup B_2$, 
so that the amalgamation of~$B_1$ and~$B_2$ over~$A$ can be taken to 
be the set of ideals of  $B_1\cup B_2$. 
(See \cite[Section V.4]{Gratzer:1998} for a detailed version of this proof.)

% We will not need a ``canonical'' amalgamation; any extension will do.
% However, our lattices $B_1$ and $B_2$ will usually be finite, in which case
% we  also expect $D$ to be finite.  

\begin{figure}[h]
\begin{picture}(0,0)%
\includegraphics{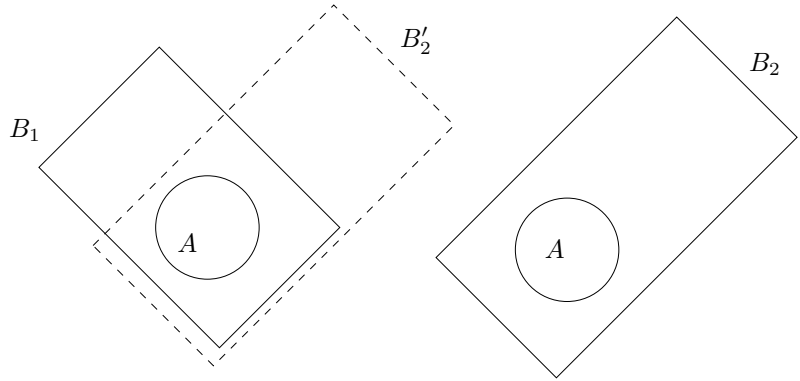}%
\end{picture}%
\setlength{\unitlength}{3315sp}%
\begingroup\makeatletter\ifx\SetFigFont\undefined%
\gdef\SetFigFont#1#2#3#4#5{%
  \reset@font\fontsize{#1}{#2pt}%
  \fontfamily{#3}\fontseries{#4}\fontshape{#5}%
  \selectfont}%
\fi\endgroup%
\begin{picture}(5922,2814)(2461,-5698)
\put(2476,-3886){\makebox(0,0)[lb]{\smash{{\SetFigFont{10}{12.0}{\rmdefault}{\mddefault}{\updefault}{$B_1$}%
}}}}
\put(3736,-4741){\makebox(0,0)[lb]{\smash{{\SetFigFont{10}{12.0}{\rmdefault}{\mddefault}{\updefault}{$A$}%
}}}}
\put(5401,-3211){\makebox(0,0)[lb]{\smash{{\SetFigFont{10}{12.0}{\rmdefault}{\mddefault}{\updefault}{$B'_2$}%
}}}}
\put(6481,-4786){\makebox(0,0)[lb]{\smash{{\SetFigFont{10}{12.0}{\rmdefault}{\mddefault}{\updefault}{$A$}%
}}}}
\put(8011,-3391){\makebox(0,0)[lb]{\smash{{\SetFigFont{10}{12.0}{\rmdefault}{\mddefault}{\updefault}{$B_2$}%
}}}}
\end{picture}%

\caption{Amalgamation with inclusion}\label{abcd}
\end{figure}

Combining this fact with the finite embeddability property we easily see: 
\begin{Proposition} The amalgamation of two finite lattices is (or: can be 
chosen to be) finite.    In other words:  The class of finite
lattices has the amalgamation property. 
\end{Proposition}
If $f_1$ and $f_2$ are the identity embeddings, then we can 
also identify $B_1$ with $g_1(f_1(B_1)) = g_1(B_1)$ and get the following 
convenient reformulation: 
\begin{Proposition}\label{ABCD}
Let $\LL$ be a class of algebras with the amalgamation property which 
is closed under isomorphism. Let
$A\le B_1$, $A\le B_2$, and $A,B_1, B_2$ be algebras in~$\LL$.  

Then we can find algebras  $B'_2$ and $D$ in~$\LL$ with (see Figure~\ref{abcd}): 
$$ B'_2\le D, \  A \le B_1 \le D , \ \mbox{ and } B_2\simeq_A B'_2.$$
\end{Proposition}

\begin{Definition}\label{JEP}
We say that a family $\LL$ of structures has the joint embedding property if
for any structures $A,B\in \LL$ there is a structure $C\in \LL$ which contains 
isomorphic copies of both $A$ and $B$.  If $\LL$ contains the empty 
structure (or more generally: has an initial element),
then the amalgamation property of~$\LL$ implies the joint embedding
property for~$\LL$.
\end{Definition}
\begin{Definition}\label{F.class}
We call a  family $\LL$ of structures  a 
\emph{Fra{\"\i}ss\'e class} if $\LL$ is closed under substructures, and has the 
amalgamation property and the joint embedding property.
\end{Definition}

\begin{Proposition} The family of lattices and the family of finite 
lattices are Fra{\"\i}ss\'e classes. 
\end{Proposition}

\subsubsection{Ultrahomogeneity and the Fra{\"\i}ss{\'e} limit}

\begin{Definition}\label{U12}
Let $\K$ be a class of finite algebras, and let $L$ be a countable algebra.
We say that $L$ is a \emph{Fra{\"\i}ss\'e structure} for~$\K$ if
\\[-3mm]
\begin{itemize}
\item [(U1)]
\begin{minipage}[t]{5cm}
For all $A,B\in \K$:\\  If $A\le B$ and $A \le L$, 
\\
then there is $B'\in \K$ with $A \le B' \le L$
and $B\simeq_A B'$. 
% See Figure~\ref{ultra}.
\\ (We call this property ``ultrahomogeneity''; note that
this word is often used for the  following related notion, see
Theorem~\ref{F.limit}(b): 
Every finite partial isomorphism between [partial] substructures
of~$L$ extends to an automorphism of~$L$.)  
\end{minipage} 
\quad
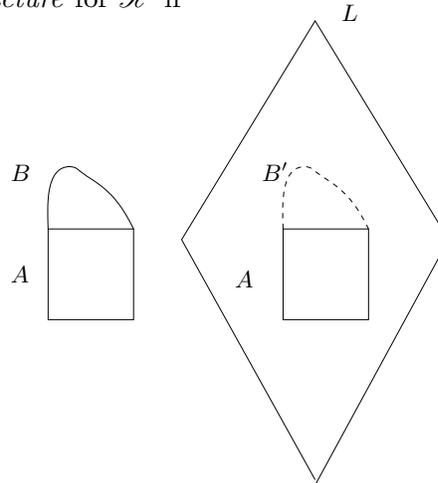
\begin{minipage}[t]{8cm}
\vskip-9mm
\vtop{\hsize=8cm \hrule height0cm depth 0cm
\begin{picture}(0,0)%
\includegraphics{ultra.pstex}%
\end{picture}%
\setlength{\unitlength}{2942sp}%
\begingroup\makeatletter\ifx\SetFigFont\undefined%
\gdef\SetFigFont#1#2#3#4#5{%
  \reset@font\fontsize{#1}{#2pt}%
  \fontfamily{#3}\fontseries{#4}\fontshape{#5}%
  \selectfont}%
\fi\endgroup%
\begin{picture}(3717,4086)(2146,-5338)
\put(4051,-3661){\makebox(0,0)[lb]{\smash{{\SetFigFont{9}{10.8}{\rmdefault}{\mddefault}{\updefault}{$A$}%
}}}}
\put(4951,-1411){\makebox(0,0)[lb]{\smash{{\SetFigFont{9}{10.8}{\rmdefault}{\mddefault}{\updefault}{$L$}%
}}}}
\put(4276,-2761){\makebox(0,0)[lb]{\smash{{\SetFigFont{9}{10.8}{\rmdefault}{\mddefault}{\updefault}{$B'$}%
}}}}
\put(2161,-3616){\makebox(0,0)[lb]{\smash{{\SetFigFont{9}{10.8}{\rmdefault}{\mddefault}{\updefault}{$A$}%
}}}}
\put(2161,-2761){\makebox(0,0)[lb]{\smash{{\SetFigFont{9}{10.8}{\rmdefault}{\mddefault}{\updefault}{$B$}%
}}}}
\end{picture}%
\hrule height0cm depth 0cm
}
\end{minipage}
\item [(U2)] $L$ contains isomorphic copies of all
elements of~$\K$.\\ (This is sometimes written as ``$L$ is universal''. 
If $\K$ contains the empty algebra or more generally has an
initial element, then this property follows from U1.)
\item [(U3)] All finite substructures of~$L$  are in~$\K$. 
\end{itemize}
\end{Definition}

If $\K$ is a  Fra{\"\i}ss\'e class of finite algebras,
then we can inductively build a countable Fra{\"\i}ss\'e structure for~$\K$.  
We construct 
an increasing sequence $(C_n:n =1,2,\ldots)$ of finite algebras such that 
\begin{itemize}
\item for all $A,B\in \K$ and all $n$ we have:  if $A\le B$ and $A\le C_n$,
      then there is $n'>n$ and $A\le B'\le C_{n'}$ with $B \simeq_A B'$; 
\item for all $B\in \K$ there is  $n'$ and $ B'\le C_{n'}$
                with $B \simeq B'$. 
\end{itemize}
The result of such a construction is called the
\emph{Fra{\"\i}ss{\'e} limit} of
$\K$.  It is unique up to isomorphism:

\begin{theorem}[{\cite{Fraisse:2000},   
see also \cite[Chapter 6.1]{Hodges:1997}}]   
\label{F.limit} 
Let $\K$ be a Fra{\"\i}ss\'e  class of finite algebras.
Then:
\begin{itemize}
\item[(a)]
 There is a unique locally finite countable structure $L^*$ which 
is a Fra{\"\i}ss\'e structure for~$\K$.
\item[(b)]
Moreover, if $L_1$ and $L_2$ are both locally finite Fra{\"\i}ss\'e
structures 
for~$\K$, then every partial isomorphism between finite subalgebras
of~$L_1$ and~$L_2$ can be extended to an isomorphism from $L_1$
onto $L_2$.
\item[(c)] Every countable locally finite algebra is isomorphic to 
a subalgebra of~$L^*$, assuming that all its finite subalgebras are in 
$\K$.
\end{itemize}
\end{theorem}

% \newpage
\subsection{Topology}

Let $X$ be a complete metric space.  A subset $M \subseteq X$ is
called ``nowhere dense'' if there is no open set contained in the
topological closure of~$M$, and $M$ is called meager (or: ``of first
category''), if $M$ can be covered by countably many nowhere dense
(or: nowhere dense closed) sets. 

Clearly, the family of meager sets forms an ideal. Since $X$ is a
complete metric space, Baire's theorem tells us that this ideal is
proper:  $X$ itself is not meager. (In fact, no nonempty open set is
meager.) 

\subsubsection{A metric space of algebras}\label{metric.space}

Recall that for notational simplicity we only consider
algebras with 2 binary operations.

We now consider all algebras (with two binary operations) whose underlying 
set is the set $\N$ of natural numbers.    A binary operation on~$\N$
is just a map from $\N\times \N$ into $\N$, i.e., an element
of~$\N^{\N^2}$.  
Identifying each algebra $A = (\N, \vee_A, \wedge_A)$ with the pair 
$(\vee_A, \wedge_A)\in \N^ {\N^2} \times \N^ {\N^2} $, we see that 
our set of algebras is really the set
$\N^ {\N^2} \times \N^ {\N^2} $, which we abbreviate as $\NN$.

The space $\NN$ in the next proposition is
sometimes called ``Baire space''; it
is homeomorphic to the set of irrational real numbers.

\begin{Proposition}\label{topology} Using the discrete topology on~$\N$, and
the Tychonoff topology (product topology) on 
$\N^ {\N^2} $ and  
$\NN=\N^ {\N^2} \times \N^ {\N^2} $, the space 
$\NN$ is a perfect Polish space (i.e., separable and completely
metrizable).
%    Note  that for all $a,b,c\in \N$ the set $\{(f,g): f(a,b)=c\}$ is clopen.

An example of a complete metric on~$\NN$ is given 
by $d((f,g), (f',g')) = 2^{-n}$, where $n$ is minimal such that 
there exist $i,j\le n$ with $f(i,j)\not=f'(i,j)$ or $g(i,j)\not=g'(i,j)$.
\end{Proposition}

\begin{Proposition}\label{subspace}
Let $\LL$ be a variety. Then $\LL\cap \NN$ is a closed subset
of~$\NN$, hence also a Polish space.
\end{Proposition}

\begin{Definition}\
\begin{itemize}
\item 
$\NN:= \N^{\N^2}\times \N^{\N^2}$   is the set of all algebras on the 
 fixed countable set~$\N$. 
\item 
$\LLn:= \LL\cap \NN$. 
\item We write $\LL\lf$ for the locally finite algebras in~$\LLn$.
\item 
We write $\Lfin$ for the class of all finite algebras whose universe 
is a subset of~$\N$. 
\item We write $\tK$ for the class of all finite partial algebras which 
are weak subalgebras of an algebra in~$\LLn$. 
\end{itemize}
\end{Definition}

\begin{Definition}  ``Almost all countable algebras in~$\LL$ have property $X$''
will mean: 
\begin{quote}
The set $\{ L\in \LLn:  \mbox{$L$ does not have property $X$}\}$ is meager
\\(in the complete metric space $\LLn$). 
\end{quote}
\end{Definition}

\begin{Definition} 
For any finite partial algebra $P\in \tK$ we let 
$$ [P]_w:= \{L \in \LLn:  P\le_w  L\}. $$ 
% , \qquad [P]_r:= \{L \in \LL:  P\le_r  L\}.$$

If $P$ is a total algebra, we write $[P]$ instead of~$[P]_w$.
\end{Definition} 

Note that for
each $a,b,c\in \|P|$ the sets
$\{L\in \LLn: a \wedge_L b = c \}$ and 
$\{L\in \LLn: a \vee_L b = c \}$ are clopen, by the definition of the 
product topology in $\LLn$; hence $[P]_w $ is clopen, as a finite 
intersection of clopen sets.   A closer inspection of the open sets
in the product topology yields the following proposition:

\begin{Proposition}\label{BASIS} The family 
 $\{[P]_w: P \in  \tK\}$ is a clopen basis for the topology on~$\LLn$. 
In other words: 
For all $L\in \LLn$, every open neighborhood of~$L$ contains a neighborhood 
of the form $[P]_w$, for some $P\in \tK$, $P\le_w L$.
(Here, $P$ is a finite \emph{partial} algebra.)

For locally finite $L$ we can ignore partial algebras altogether: 
If $L\in \LL\lf$, then every open neighborhood of~$L$ contains a neighborhood 
of the form $[A]$, for some $A\in \Lfin$, $A\le L$. 
(Here, $A$ is a finite \emph{total} algebra.)
\end{Proposition}

\section{The typical countable algebra}\label{proofs}

Throughout this section the class $\LL$ will be the class of countable algebras
of some fixed variety.

We write 
$\LLn$ for the set of all algebras in~$\LL$ whose
carrier set is the set $\N$ of  natural numbers.
At various points we will 
assume that $\LL$ has several properties, in particular:
the finite embeddability property, the amalgamation property and the joint embedding property.  If we take $\LL$ to be the class of countable lattices, then 
$\LL$  has all these properties. 

\subsection{Two residual subsets of~$\LLn$}
Our first lemma shows that the f.e.p.\ allows us to  restrict our attention to 
locally finite algebras; thus the neighborhoods in our topological space
are generated already by the finite total algebras and we may ignore
partial algebras.

\begin{Lemma}\label{lf}
Assume that $\LL$ has the 
finite embeddability property.
Then
almost all algebras $L\in \LLn$ are locally finite.  

In fact, the set $\LL\lf$ of locally finite algebras  is the intersection 
of countably many dense open sets of~$\LLn$ (and hence, with the
subspace topology, also a Polish space).
\end{Lemma}
\begin{proof}
We have $\displaystyle\LL\lf = \bigcap\limits_{k\in \N}
\bigcup\limits_{
\atopp
{\{0,\ldots, k\}\subseteq \|B|}{  B \in \Lfin}} [B]$. 

We claim that each set 
$$\L_k:= \{ L \in \LLn : \exists B\in \Lfin,\,  B\le L, \, \{0,\ldots, k\} 
\subseteq \|B|\}$$
 is open dense.  It is clear that $\L_k$ is open. 

Now fix any  nonempty open subset of~$\LLn$.   By Proposition~\ref{BASIS}
we may assume that this open set is of the form $[P]_w$, where $P\in\tK$ is a
partial algebra. Find a partial algebra $P'\in\tK$, $P\le_w P'$
such that $\|P'|$
contains the set $\{0,\ldots, k\}$.   By the finite embeddability property,
there is a total algebra  $B\in \Lfin$, 
$P' \le B$.  So $[P]_w \supseteq [P']_w \supseteq [B]$.  $[B]$ is 
nonempty,  so $\L_k$
meets $[P]_w$.   So $\L_k$ is dense. 
\end{proof}

\begin{Proposition}
Assume that the class $\LL$ is a Fra{\"\i}ss\'e class with the 
finite embeddability property.
Then the 
 set 
 $$\U:=\{L\in \LL\lf: L \mbox{ is Fra{\"\i}ss\'e for~$\LL$}\}$$ is
residual in~$\LL\lf$, hence also in~$\LLn$.
\end{Proposition}
\begin{proof} We can write $\U$ as $\U=\U_1\cap \U_2$, where 
for~$\ell=1,2$ the set $\U_\ell$ is the set of all $L\in \LLn$ satisfying 
property U$\ell$ in Definition~\ref{U12}. We will show that $\U_1$ is residual,
leaving the case of~$\U_2$ to the reader. 

We have 
$$\begin{aligned} \U_1&= \bigcap_{A\le B\in \Lfin}  \{L : 
\text{ if } A\le L \text{ then } \exists B'\le L: B \simeq_A B'\} = \\
&=
\bigcap_{A\le B\in \Lfin} (\U_{\lnot A} \cup \U_{A,B}),
\end{aligned}
$$
where 
$ \U_{\lnot A}=\LL\lf\setminus [A]$
is the set of all $L\in \LL\lf$ with $A\not\le L$, 
and $$\U_{A,B}:= \{L\in \LL\lf: A \le L, \mbox{ and } 
\exists B'\le L:  B\simeq_A B'\,\}.$$
The sets $\U_{\lnot A}$ and  $\U_{A,B}$  are open ($\U_{A,B}$ is 
the union of sets of the form $[B']$, $B'\in \Lfin$).

We now check that each set $ \U_{\lnot A} \cup \U_{A,B}$ is dense. 

So let $[C] $ be a basic open set in~$\LL\lf$, $C\in \Lfin$.  We may
assume  that $\|A|\subseteq \|C|$.   If $A\not\le C$, then  $[C]
\subseteq \U_{\lnot A}$, so we are done. 

So assume $A\le C$.   As also $A\le B$, we can find a finite algebra 
$D$ and a subalgebra $B'$ such that $B \simeq_A B'$ (using
Proposition~\ref{ABCD}).
 Wlog $\|D|\subseteq \N$.   We have $[D] \subseteq
[C]$, 
and $[D] \subseteq \U_{A,B}$. 
\end{proof}
\begin{Corollary}
Assume that $\LL$ is a Fra{\"\i}ss\'e class with 
the joint embedding property.  Let $L^*$ be the Fra{\"\i}ss\'e 
limit of the finite algebras in~$\LL$.  Then almost all  countable 
algebras in~$\LLn$ are Fra{\"\i}ss\'e structures for the finite structures
in~$\LL$,
and hence are isomorphic to~$L^*$.
\end{Corollary}

In particular, taking $\LL=$ the class of all 
countable lattices, we get: There is 
a (``typical'' or ``generic'')
countable lattice $L^*$ such that almost all countable lattices
are isomorphic to~$L^*$.

% \newpage
\section{Examples} \label{examples}
Assume that $\LL$ has the finite embeddability property, the amalgamation
property and the joint embedding property.  Let $L^*$ be the Fra{\"\i}ss\'e 
limit of the finite algebras in~$\LL$.

We have already seen that $L^*$ is locally finite and contains copies 
of all finite $\LL$-algebras, and (by Theorem~\ref{F.limit}(c)) 
even all locally finite countable 
$\LL$-algebras.
  The finite embeddability property implies 
that $L^*$ satisfies no law that does not hold in all $\LL$-algebras; 
in other words, $L^*$ generates $\LL$.

Since every finite automorphism between subalgebras (even partial subalgebras)
extends to an automorphism of~$L^*$, $L^*$ has a very rich automorphism 
group. 

The typical countable Boolean algebra is the (unique)
atomless countable Boolean algebra;  the typical distributive lattice 
and its automorphism group have been investigated
 in \cite{Droste+Macpherson:2000}. 

We now consider two special cases which are not locally finite: lattices
and bounded lattices. 

\subsection{Lattices}

Let $L^*$ be the typical countable lattice. 

As remarked in the introduction, $L^*$ is certainly not distributive (it
contains~$N_5$), and in fact satisfies no law that is not implied by laws
defining the variety of lattices.     Since $L^*$ contains all countable
locally finite lattices, $L^*$ contains, for example, a chain isomorphic to 
the rational numbers, and also an infinite antichain.

Note that $L^*$ is very different from the random order $R$, i.e., 
the Fra{\"\i}ss\'e limit of the class of finite
partial orders: in $R$ there are
\emph{no} no elements  $x \not=y$ which have  a smallest upper bound.

\begin{Proposition} For every monotone function~$f:(L^*)^n\to L^*$ and every finite
sublattice $A\le L^*$ there is a lattice polynomial $p(x_1,\ldots, x_n)$
which induces the function~$f $ on~$A$.   (We say that $L^*$ has the 
\emph{monotone interpolation property}.)

This implies that $L^*$ is simple (i.e., has no nontrivial congruence
relations). 
\end{Proposition}
\begin{proof}
By \cite{mfl}, we can find a finite lattice $B$, $A\le B$ and a polynomial
$p$ with 
coefficients in~$B$ which interpolates $f$ on~$\|A|$, i.e., $p(\vec a) = 
f(\vec a)$ for all $\vec a\in A^n$. 
We can use the ultrahomogeneity of~$L^*$ to find~$B'\le L^*$,
$B\simeq _A B'$.  The isomorphism between~$B$ and~$B'$ translates
$p$ to a polynomial~$p'$ with coefficients in~$B'$
which still interpolates~$f$ on~$A$.

The fact that the monotone interpolation property 
implies that there are no nontrivial
congruence relations is well known: for any nontrivial congruence relation 
$\sim$ we can find $a_1 < a_2$ and $b_1 <b_2$
such that $a_1\sim a_2$, $b_1\not\sim b_2$. There 
is a monotone total function mapping $a_i$ to $b_i$; as such a function 
does not respect $\sim$,  such a 
function cannot be a polynomial. 
\end{proof}

\begin{Proposition} 
Any two nontrivial intervals in the typical lattice are isomorphic to
each other. (In fact, the isomorphism can be taken to be the restriction
of an automorphism of~$L^*$.) 
\end{Proposition}
\begin{proof} By  Theorem~\ref{F.limit}(b). 
\end{proof}

\subsection{$\zo$-lattices}

We now consider a language where we have  two constant symbols $0,1$ 
in addition to the two operations $\wedge$ and $\vee$.   We will consider
the variety of~$\zo$-lattices (i.e., lattices in which $0$ and $1$ 
are the greatest and smallest element). 

Note that the 1-element $\zo$-lattice does not embed into any other
$\zo$-lattice, so we will have to ignore it in our considerations.
{\bf From now on, \{0,\,1\}-lattice will mean: \{0,\,1\}-lattice with
$\bf 0\not=1$.}
 
 The same construction as in Proposition~\ref{lat.amalg} 
 (after the necessary change of notation, taking into account the new 
 constants) shows the following: 
 \begin{Proposition}\label{01.amalg}
 Whenever $A$ is a $\zo$-lattice, 
 and $f_1:A\to B_1$, $f_2:A\to B_2$ are $\zo$-homomorphisms, then there exists a $\zo$-lattice $D$ and $\zo$-homomorphisms $g_i:B_i\to D$ such that 
 $g_1 \circ f_1 = g_2 \circ f_2$.  Moreover, if $B_1$ and $B_2$ are
 finite, then $D$ can be chosen to be finite. 

 In other words:  the class of $\zo$-lattices, as well as the class
 of  finite $\zo$-lattices, has the amalgamation property.
 \end{Proposition}
 \begin{Proposition} The class of~$\{0,1\}$-lattices has the finite embeddability property. 
 \end{Proposition}
 \begin{proof} Let $A $ be a partial finite $\{0,1\}$-lattice.  Wlog
 we may assume that $0$ and $1$ are defined in $A$.  The natural
 embedding  
 of~$A$ into the lattice of \emph{nonempty} ideals preserves $0$ and
 $1$ (since $0$ is mapped to the smallest nonempty ideal~$\{0\}$).
 \end{proof}

It is now easy 
to see that the  $\zo$-lattices are a Fra{\"\i}ss\'e class with the f.e.p. 
Hence there is a ``typical'' countable $\zo$-lattice $K^*$. 
Note that $1$ is typically not join-irreducible (i.e., there are lattices,
even finite ones, where $1 = x\vee y$ for some $x,y<1$), so $1$ is also
not join-irreducible in $K^*$.  

Let $L^*$ be the typical lattice, and let $a,b\in L^*$, $a<b$. Consider the 
interval $[a,b]$ as a $\{0,1\}$-lattice with $a=0$, $b=1$. 

The ultrahomogeneity/universality
of~$L^*$ (with respect to the class of finite lattices)
easily implies that $[a,b]$ is ultrahomogeneous/universal (with respect 
to the class of finite $\zo$-lattices).  Hence  we get: 

\begin{Proposition}
Every nontrivial interval in the typical lattice $L^*$ is isomorphic 
to every nontrivial interval in the typical $\zo$-lattice $K^*$, in particular 
to $K^*$ itself.
\end{Proposition}

Using \cite{l01} instead of \cite{mfl} we can also show:

\begin{Proposition} For every monotone function~$f:(K^*)^n\to K^*$ and every finite
$\zo$-sublattice $A\le K^*$ there is a lattice polynomial $p(x_1,\ldots, x_n)$
which induces the function~$f $ on~$A$. 

This again implies that $K^*$ is simple. 
\end{Proposition}

\section{Acknowledgements}
I am grateful to Manfred Droste for pointing out several papers 
on ubiquitous structures, to Michael Pinsker for a conversation that led to 
Lemma~\ref{lem.equiv}, to Günther Eigenthaler for his careful proofreading,
and to both referees for their constructive
criticism. 

% \newpage
\bibliography{goldstrn,other} \bibliographystyle{plain}

\bigskip\bigskip \def\nbibitem[#1]{\advance\bcount1 \bibitem[\number\bcount]}
  \newcount\bcount \def\acceptedinprint{accepted/in print}
  \def\submittedinpreparation{submitted/in preparation}
\begin{thebibliography}{Goldstern 1996}

\bibitem[Bankston-Ruitenburg]{Bankston+Ruitenburg:1990}
Paul Bankston and Wim Ruitenburg.
\newblock Notions of relative ubiquity for invariant sets of relational
  structures.
\newblock {\em J. Symbolic Logic}, 55(3):948--986, 1990.

\bibitem[Baumslag-Solitar]{Baumslag-Solitar:1962}
Gilbert Baumslag and Donald Solitar.
\newblock Some two-generator one-relator non-{H}opfian groups.
\newblock {\em Bull. Amer. Math. Soc.}, 68:199--201, 1962.

\bibitem[Droste-Kuske]{Droste+Kuske:2007}
Manfred Droste and Dietrich Kuske.
\newblock Almost every domain is universal.
\newblock {\em Electronic Notes in Theoretical Computer Science}, 173:103--119,
  2007.


\bibitem[Droste-Macpherson]{Droste+Macpherson:2000}
Manfred Droste and Dugald Macpherson.
\newblock The automorphism group of the universal distributive lattice.
\newblock {\em Algebra Universalis}, 43(4):295--306, 2000.


\bibitem[{Fra{\"\i}ss\'e}]{Fraisse:2000}
Roland Fra{\"{\i}}ss{\'e}.
\newblock {\em Theory of relations}, volume 145 of {\em Studies in Logic and
  the Foundations of Mathematics}.
\newblock North-Holland Publishing Co., Amsterdam, revised edition, 2000.
\newblock With an appendix by Norbert Sauer.

\bibitem[Funayama]{Funayama:1953}
Nenosuke Funayama.
\newblock Notes on lattice theory. {IV}. {O}n partial (semi-) lattices.
\newblock {\em Bull. Yamagata Univ. (Nat. Sci.)}, 2:171--184, 1953.

\bibitem[Goldstern 1996]{mfl}
Martin Goldstern.
\newblock {Interpolation of monotone functions in lattices}.
\newblock {\em Algebra Universalis}, 36:108--121, 1996.

\bibitem[Goldstern 1997]{most}
Martin Goldstern.
\newblock {Most algebras have the interpolation property}.
\newblock {\em Algebra Universalis}, 38:97--114, 1997.

\bibitem[Goldstern 1998]{l01}
Martin Goldstern.
\newblock Interpolation of monotone functions in $\{0,1\}$-lattices.
\newblock In {\em Contributions to General Algebra 10}. Heyn Verlag, 1998.

\bibitem[Gr\"atzer 1979]{Gratzer:1979}
George Gr{\"a}tzer.
\newblock {\em Universal algebra}.
\newblock Springer-Verlag, 1979.

\bibitem[Gr\"atzer 1998]{Gratzer:1998}
George Gr{\"a}tzer.
\newblock {\em {General lattice theory. 2nd ed.}}
\newblock Birkh\"auser, Basel, 1998.

\bibitem[Hodges]{Hodges:1997}
Wilfrid Hodges.
\newblock {\em A shorter model theory}.
\newblock Cambridge University Press, 1997.

\bibitem[J\'onsson]{Jonsson:1956}
Bjarni J{\'o}nsson.
\newblock Universal relational systems.
\newblock {\em Math. Scand.}, 4:193--208, 1956.

\bibitem[Lindner-Evans]{Lindner-Evans:1977}
Charles~C. Lindner and Trevor Evans.
\newblock {\em Finite embedding theorems for partial designs and algebras}.
\newblock Les Presses de l'Universit\'e de Montr\'eal, Montreal, Que., 1977.
\newblock S\'eminaire de Math\'ematiques Sup\'erieures, No. 56 (\'Et\'e 1971).

\bibitem[Robinson]{Robinson:1996}
Derek J.~S. Robinson.
\newblock {\em A course in the theory of groups}, volume~80 of {\em Graduate
  Texts in Mathematics}.
\newblock Springer-Verlag, New York, second edition, 1996.

\end{thebibliography}

\end{document}